\documentclass[12pt]{amsart}

\usepackage{epsfig}

\usepackage{amssymb}
\usepackage{amsmath}

\def\ds{\displaystyle}
\def\bs{\boldsymbol}

\newcommand{\R}{{\mathbb{R}}}

\newcommand{\SF}{{\mathbb{S}}}

\def\tl{\widetilde}

\def\bu{{\mathbf{u}}}
\def\bv{{\mathbf{v}}}

\def\bff{{\mathbf{f}}}

\def\A{{\mathbf{A}}}

\def\CM{\bs{\mathfrak{C}}}

\def\K{\bs{\mathfrak{K}}}

\def\LLb{{\bs{\mathfrak{L}}}}

\def\l{\ell}
\def\LL{{{\mathfrak{L}}}}

\newtheorem{thm}{Theorem}[section]
 \newtheorem{cor}[thm]{Corollary}
 \newtheorem{lem}[thm]{Lemma}
 \newtheorem{prop}[thm]{Proposition}
 \theoremstyle{definition}
 \newtheorem{defn}[thm]{Definition}
 \theoremstyle{remark}

 \numberwithin{equation}{section}

\begin{document}

\title[Discontinuous Elliptic Systems]
 {Fine Regularity for Elliptic Systems\\ with Discontinuous Ingredients}

\author[D. Palagachev]{Dian  Palagachev }

\address{
   Politecnico di Bari\\
   Dipartimento di Matematica\\
   Via E. Orabona, 4\\
   70~125~Bari,\ Italy}

\email{dian@dm.uniba.it}

\author[L. Softova]{Lubomira  Softova}
\address{Bulgarian Academy of Sciences\\
   Institute of Mathematics\\
   Sofia,\ Bulgaria}
\email{luba@dm.uniba.it}

\subjclass{Primary 35J45; Secondary 35R05, 35B45, 35B65, 46E35}

 \keywords{Elliptic systems, a'priori estimates, H\"older regularity, singular
integrals}


\begin{abstract}
We  propose results on interior Morrey, $BMO$ and H\"older regularity
 for the strong solutions to linear elliptic systems of order $2b$ with
discontinuous coefficients
and right-hand sides belonging to the Morrey space $L^{p,\lambda}.$
\end{abstract}

\maketitle

\section{ Main Results}

It is well known that,
 when dealing with elliptic systems, in contrast to the case of a single second-order
elliptic equation,
 the solely essential boundedness of the principal coefficients is not
sufficient to ensure
 H\"older continuity even of the solution (see \cite[Chapter~1]{MPS}). On the
other hand,
 precise estimates on  H\"older's seminorms of the solution and its lower order
derivatives is a matter of great
 concern in the study of nonlinear elliptic systems. In fact, these bounds imply
good
 mapping properties of certain Carath\'eodory operators ensuring this way the
possibility to
 apply the powerful tools of the nonlinear analysis and differential calculus.
It turns out that
 ``suitable continuity'' of the principal coefficients of the system under
consideration is sufficient
 to guarantee good regularity (e.g. Sobolev) of the solutions (see
\cite{DN,Sl}). We deal here with
 discontinuous coefficients systems for which the discontinuity is expressed in
terms of appurtenance
 to the class of functions with vanishing mean oscillation. Although such
systems have been already
 studied in Sobolev spaces $W^{2b,p}$ (cf. \cite{CFF}), our functional framework
 is that of the Sobolev--Morrey classes $W^{2b,p,\lambda}.$ These possess better
embedding properties into H\"older spaces
 than $W^{2b,p}$ and as outgrowth of suitable Caccioppoli-type estimates, we
give precise characterization
of the Morrey, $BMO$ or H\"older regularity of the solution
and its derivatives up to order $2b-1.$

Let $\Omega$ be a
domain in $\R^n,$ $n\geq2,$ and  consider the  linear system
\begin{equation}\label{system}
 \LLb (x,D) \bu:= \sum_{|\alpha|=2b}\A_\alpha (x) D^\alpha \bu(x) =
 \bff(x)
\end{equation}
 for the unknown vector-valued function $\bu\colon \Omega\to\R^m$ given by the
transpose $\bu(x)=\big(u_1(x),\ldots,u_m(x)\big)^{\mathrm{T}},$ $m\geq 1,$
$\bff(x)=(f_1(x),\ldots,f_m(x))^{\mathrm{T}},$ where
 $\A_\alpha(x)$ is the $m\times m$-matrix $\left\{ a_\alpha^{jk}(x)
\right\}_{j,k=1}^m$
and $a_\alpha^{jk}\colon \Omega\to\R$ are measurable functions.
 Hereafter, $b\geq 1$ is a fixed integer, $\alpha=(\alpha_1,\ldots,\alpha_n)$ is
a multiindex of length $|\alpha|=\alpha_1+\cdots+\alpha_n$ and
$D^\alpha:=D_1^{\alpha_1}\ldots
D_n^{\alpha_n}$ with $D_i:=\partial/\partial x_i.$ This way, the matrix
differential operator
 $\LLb (x,D)$ has entries 
$$\l^{jk}(x,D):= \sum_{|\alpha|=2b}a_\alpha^{jk}(x)D^\alpha
$$
 and for fixed
 $j$ and $k$ the polynomial
 $$
\l^{jk}(x,\xi):=\sum_{|\alpha|=2b}a_\alpha^{jk}(x)\xi^\alpha,\quad \xi\in\R^n,\quad
\xi^\alpha :=\xi_1^{\alpha_1}\xi_2^{\alpha_2}\cdots\xi_n^{\alpha_n},
$$
is homogeneous of degree $2b.$

We suppose \eqref{system} to be an
{\it elliptic system,\/} that is, the characteristic determinant of
$\LLb (x,\xi)$
 is non-vanishing for a.a. $x\in\Omega$ and all $\xi\neq0.$ In view of the
homogeneity of $\l^{jk}$'s, this rewrites as (see  \cite{CFF,DN})
\begin{equation}\label{elliptic}
 \exists\  \delta>0\colon\quad \text{det\,} \Big\{\sum_{|\alpha|=2b}
\A_\alpha(x)\xi^\alpha
 \Big\}\geq \delta |\xi|^{2bm}\quad \text{a.a. }x\in\Omega,\  \forall \xi \in
\R^n.
\end{equation}

Our goal  is to obtain interior  H\"older regularity of the solutions to
\eqref{system}
as a  byproduct of a'priori estimates in Sobolev and Sobolev-Morrey spaces. Let
us recall the definitions of these functional classes.

\begin{defn}\label{dSob}\em
The  Sobolev space $W^{2b,p}(\Omega),$ $p\in (1,+\infty),$  is the
collection of $L^p(\Omega)$ functions $u\colon \Omega\to\R$ all of which
distribution derivatives  $D^\alpha u$ with $|\alpha|\leq 2b,$
belong to
$L^p(\Omega).$ The norm in $W^{2b,p}(\Omega)$ is
$$
 \|u\|_{W^{2b,p}(\Omega)}:=\ \sum_{s=0}^{2b} \sum_{|\alpha|=s} \|D^\alpha u\|_{p;\Omega},\qquad
 \|\cdot\|_{p;\Omega}:=\ \left(\int_\Omega|\cdot|^pdx \right)^{1/p}.
$$
 For the sake of brevity, the cross-product of $m$ copies of $L^p(\Omega)$ is
denoted by the same symbol. Thus, if
$\bu=(u_1,\ldots,u_m)$ is a vector-valued function,
 $\bu\in L^{p}(\Omega)$ means that $u_k\in L^p(\Omega)$ for all
$k=1,\ldots,m,$ and
$\|\mathbf{u}\|_{p;\Omega}:=\sum_{k=1}^m \|{u_k}\|_{p;\Omega}.$
\end{defn}

\begin{defn}\label{dMor}\em
Let $p\in(1,+\infty)$ and $\lambda\in(0,n).$
The function $u\in L^p(\Omega)$ belongs to the	Morrey space
$L^{p,\lambda}(\Omega)$ if
\[
 \|u\|_{p,\lambda;\Omega}:=\ \left(\sup_{r>0}\frac{1}{r^\lambda}\int_{B_r\cap
\Omega} |u(x)|^p dx \right)^{1/p} <\infty
\]
where $B_r$ ranges in the set of  balls with radius $r$ in $\R^n.$
The Sobolev--Morrey space $W^{2b,p,\lambda}(\Omega)$
 consists of all functions $u\in
W^{2b,p}(\Omega)$ with generalized derivatives	$D^\alpha u,$
$|\alpha|\leq 2b,$ belonging to $L^{p,\lambda}(\Omega).$ The norm in
$W^{2b,p,\lambda}(\Omega)$ is given by
$$
 \|u\|_{W^{2b,p,\lambda}(\Omega)}:=\  \sum_{s=0}^{2b}\sum_{|\alpha|=s}\|D^\alpha
u\|_{p,\lambda;\Omega}. $$
\end{defn}
We refer the reader to \cite{Cm,Cm1,M}
for various properties of the Morrey and Sobolev--Morrey spaces.
\begin{defn}\label{dBMO}
For a locally integrable function $f\colon\ \R^{n}\to\R$ define
\[
\eta_f(R):=\sup_{r\leq R} \frac{1}{|B_r|}\int_{B_r}
 |f(y)-f_{B_r}|dy\quad \text{for every}\ R>0,
\]
where $B_r$ ranges over the  balls in $\R^n$ and
$f_{B_r}=\frac{1}{|B_r|}\int_{B_r} f(y) dy.$
Then:
\begin{itemize}
\itemsep=1pt
\item[$\bullet$]  $f\in BMO$ {\em (bounded mean oscillation\/}, see
John--Nirenberg~{\rm\cite{JN}}$)$ if $\|f\|_*:=\sup_R \eta_f(R)<+\infty.$
$\|f\|_*$ is a norm in $BMO$ modulo constant functions under which $BMO$
is a Banach space.
\item[$\bullet$] $f\in VMO$ {\em (vanishing mean oscillation\/}, see
Sarason~{\rm\cite{S}}$)$ if
$f\in BMO$ and $\lim_{R\downarrow 0}\eta_f(R)=0.$
The quantity $\eta_f(R)$ is referred to as  $VMO$-modulus of $f.$
\end{itemize}
The spaces $BMO(\Omega)$ and $VMO(\Omega),$ and $\|\cdot\|_{*;\Omega}$
 are defined in a similar manner
taking $B_r\cap \Omega$ instead of $B_r$ above.
\end{defn}

 As already mentioned, the desired  H\"older bounds will be derived on
the base of the following a'priori estimate in Sobolev--Morrey classes.
\begin{thm}\label{th3}
Suppose \eqref{elliptic}, $a_\alpha^{jk}\in VMO(\Omega)\cap L^\infty(\Omega),$
 $\bff\in L^{p,\lambda}_{\rm loc}(\Omega),$ $1<p<\infty,$ $0<\lambda<n,$ and let $\bu\in  W^{2b,p,\lambda}_{\rm loc}(\Omega)$ be a strong solution of
\eqref{system}.
 Then, for any $\Omega'\Subset\Omega''\Subset\Omega$ there is a constant $C$
depending on
 $n,$ $m,$ $b,$ $p,$ $\lambda,$ $\delta,$ $\|a^{jk}_\alpha\|_{\infty;\Omega},$
the $VMO$-moduli
$\eta_{a^{jk}_\alpha}$ of the coefficients {\em (cf. \cite{CFF,CFL})}
and $\mathrm{dist\,}(\Omega',\partial\Omega''),$ such that
\begin{equation}\label{9a}
\|\bu\|_{W^{2b,p, \lambda}(\Omega')}\leq
C\left(\|\bff\|_{p,\lambda;\Omega''}+ \|\bu\|_{p,\lambda;\Omega''} \right).
\end{equation}
\end{thm}
It  turns out, moreover, that the operator $\LLb$ improves the integrability of
solutions to \eqref{system}. In fact, by means of standard
homotopy arguments and making use of formula \eqref{10} (cf. \cite{CFF},
\cite[Section~3]{PS1}), it is easy to get
\begin{cor}\label{cr1}
Under the hypotheses of Theorem~$\ref{th3},$ suppose  $\bu\in W^{2b,q}_{\rm
loc}(\Omega)$
with $q\in (1,p].$ Then $\bu\in W^{2b,p,\lambda}_{\rm loc}(\Omega).$
\end{cor}

A  combination of \eqref{9a} with the embedding properties of Sobolev--Morrey
spaces leads to a
precise characterization of the Morrey, $BMO$ and H\"older regularity of the
solutions to \eqref{system}.
\begin{cor}\label{cr4}
Under the hypotheses of Theorem~$\ref{th3}$ define
$s_0$ as the least non-negative integer such that $\frac{n}{2b-s_0}>1$
and fix an $s\in\{s_0,\ldots,2b-1\}.$ Then there is a constant $C$ such that:
\begin{itemize}
\leftskip=-18pt
\item[a)] if $p\in \left(1,\frac{n-\lambda}{2b-s}\right)$ then
$D^s\bu \in L^{p,(2b-s)p+\lambda}(\Omega')$  and
$$
\|D^s\bu\|_{p,(2b-s)p+\lambda;\Omega'} \leq
C\left(\|\bff\|_{p,\lambda;\Omega''}+\|\bu\|_{p,\lambda;\Omega''}\right);
$$
\item[b)] if $p=\frac{n-\lambda}{2b-s}$ then $D^s\bu \in BMO(\Omega')$ and
$$
 \|D^s\bu\|_{BMO;\Omega'}\leq C
\left(\|\bff\|_{p,\lambda;\Omega''}+\|\bu\|_{p,\lambda;\Omega''}\right); 
$$
\item[c)] if
$p\in\left(\frac{n-\lambda}{2b-s},\frac{n-\lambda}{2b-s-
1}\right)$\footnote{This rewrites as $p\in
(n-\lambda,\infty)$ when $s=2b-1.$}  then
$D^s\bu \in C^{0,\sigma_s}(\Omega')$ with
$\sigma_s=2b-s-\frac{n-\lambda}{p}$	and
$$
\sup_{\underset{x,\,x'\in \Omega'}{x\neq x'}}
\frac{\left|D^s\bu(x)-D^s\bu(x')\right|}{|x-x'|^{\sigma_s}} \leq
C\left(\|\bff\|_{p,\lambda;\Omega''}+\|\bu\|_{p,\lambda;\Omega''}\right).
$$
\end{itemize}
  If $s_0\geq 1$ (i.e., $2b\geq n$) and $p\in
\left(1,\frac{n-\lambda}{2b-s_0}\right)$
then $\bu\in C^{s_0-1,2b-s_0+1-\frac{n-\lambda}{p}}(\Omega')$ and
$$
\sup_{\underset{x,\,x'\in \Omega'}{x\neq x'}}
  \frac{\left|D^{s_0-1}\bu(x)-D^{s_0-1}\bu(x')\right|}{|x-x'|^{2b-s_0+1-\frac{n-
\lambda}{p}}} \leq
C\left(\|\bff\|_{p,\lambda;\Omega''}+\|\bu\|_{p,\lambda;\Omega''}\right).
$$
\end{cor}

\begin{center}
\includegraphics[scale=0.60]{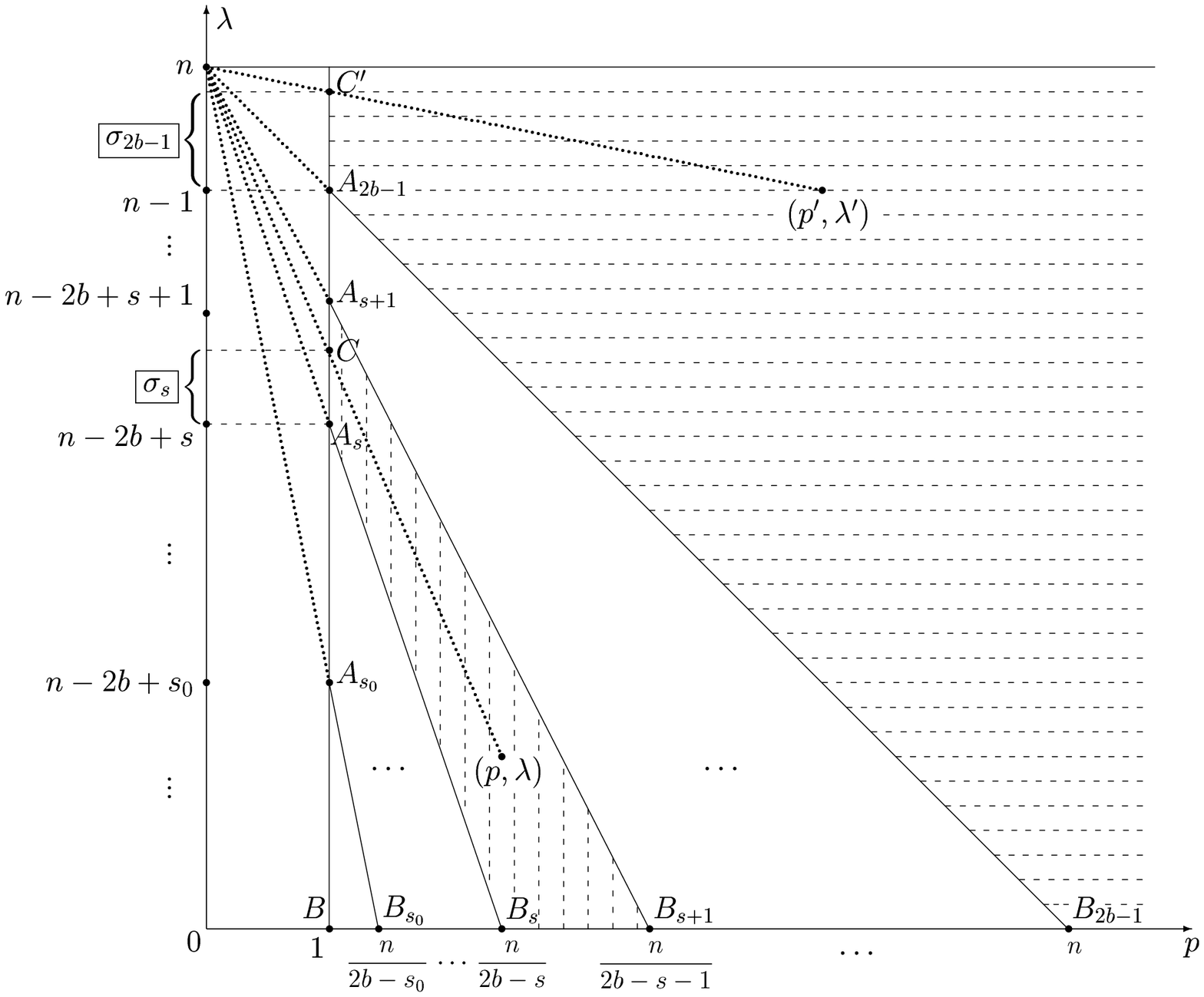}
\end{center}

The picture illustrates geometrically the results of Corollary~\ref{cr4}
with the couple $(p,\lambda)$
lying in the semistrip $\big\{(p,\lambda)\colon\ p>1,\ 0<\lambda<n\big\}$
and $s\in\{s_0,\ldots,2b-1\}.$
The points $B_s$ on the $p$-axis are $B_s=\left(\frac{n}{2b-s},0\right),$
$B=(1,0),$ and $A_s=(1,n-2b+s)$ is the intersection of
the line through $(0,n)$ and $B_s$ with the vertical line
$\{p=1\}.$

\indent
When  $(p,\lambda)$ belongs to the open right triangle $BB_sA_s$ then we have case~a), i.e.
$D^s\bu\in
L^{p,(2b-s)p+\lambda}(\Omega').$  If
$(p,\lambda)$ lies on the open segment $A_sB_s$
we have  $D^s\bu \in BMO(\Omega')$
(case~b)). In particular, $(p,\lambda)\in \triangle
BB_{s_0}A_{s_0}$ yields $D^{s_0}\bu\in L^{p,(2b-s_0)p+\lambda}(\Omega'),$
while
$\bu\in C^{s_0-1,2b-s_0+1-\frac{n-\lambda}{p}}(\Omega')$ if $s_0\geq1.$
Further  on, $(p,\lambda)\in A_{s_0}B_{s_0}$ gives $D^{s_0}\bu\in
BMO(\Omega').$

Let $s\in\{s_0,\ldots,2b-2\}$ and suppose
$(p,\lambda)$ lies in the interior of the quadrilateral
$Q_s:=B_sB_{s+1}A_{s+1}A_s.$ Then $D^s\bu\in C^{0,\sigma_s}(\Omega')$
(case~c)) whereas
$D^{s+1}\bu \in L^{p,(2b-s-1)p+\lambda}(\Omega').$
Moreover, the exponent $\sigma_s$ is the length  of the segment
${CA_s}$ where
$C=C(p,\lambda)$ is the intersection of the line $\{p=1\}$ with the
line passing through the points $(p,\lambda)$ and $(0,n).$ In particular,
$(p,\lambda)\in A_{s+1}B_{s+1}$ implies $D^{s+1}\bu\in BMO(\Omega').$

Similarly, set $Q_{2b-1}$ for the shadowed {\it unbounded region\/} on the
picture. Then $(p',\lambda')\in Q_{2b-1}$ gives $D^{2b-1}\bu\in
C^{0,\sigma_{2b-1}}(\Omega')$ with
$\sigma_{2b-1}$ equals to the length of $C'A_{2b-1},$
$C'=C'(p',\lambda'),$ while $D^{2b-1}\bu\in BMO(\Omega')$ if
$(p',\lambda')\in A_{2b-1}B_{2b-1}.$

\section{Newtonian-type Potentials}
Fix the coefficients of \eqref{system} at a point
$x_0\in\Omega$ and consider the constant coefficients operator
$\LLb(x_0,D):= \sum_{|\alpha|=2b} \A_\alpha(x_0) D^\alpha.$
The $2bm$-order differential
operator
$$L(x_0,D):= \text{det\,}\LLb(x_0,D) = \text{det\,}\Big\{\sum_{|\alpha|=2b}a_\alpha^{jk}(x_0)D^\alpha\Big\}
$$
 is
elliptic in view
of \eqref{elliptic},  and therefore there exists its fundamental solution $\tl
\Gamma(x_0;x-y)$
(see \cite{DN,J}). If the space dimension $n$ is an {\it odd\/} number, then
$$
\tl\Gamma(x_0;x-y) = |x-y|^{2bm-n}P\left(x_0;\frac{x-y}{|x-y|}\right)
$$
 where
$P(x_0;\xi)$ is a real	analytic function of
$\xi\in\SF^{n-1}:=\{\xi\in\R^n\colon |\xi|=1\}.$
In case of {\it even\/} dimension $n$ the procedure is standard and is based on
introduction of a fictitious
new variable $x_{n+1}$	and extension of all functions as constants with respect
to it, see \cite{DN} for details.
Set $\{L_{jk}(x_0,\xi)\}_{j,k=1}^m$ for the {\it cofactor matrix\/} of
$\{\l^{jk}(x_0,\xi)\}_{j,k=1}^m$  and note that for all fixed $j,k=1,\ldots,m,$ $L_{jk}(x_0,D)$ is  either a differential operator of order $2b(m-1),$ or the
operator of multiplication by $0.$
Making use of the identities 
$$
\sum_{k=1}^m
\l^{ik}(x_0,\xi)L_{jk}(x_0,\xi)=\delta_{ij}
L(x_0,\xi)
$$
 with Kronecker's $\delta_{ij},$ it is not hard to check (cf.
\cite{DN,Sl}) that the fundamental matrix $\bs{\Gamma}(x_0;x)=
\left\{\Gamma^{jk}(x_0;x)\right\}_{j,k=1}^m$ of $\LLb(x_0,D)$ has entries
$\Gamma^{jk}(x_0;x)=L_{kj}(x_0,D) \widetilde\Gamma(x_0;x).$

Let $r>0$ be so small  that $B_r=\{x\in\R^n\colon |x-x_0|<r\}\Subset\Omega,$ and
let $\bv\in C_0^\infty(B_r).$
Then, employing  
$$
\LLb(x_0,D)\bv(x)=\big(\LLb(x_0,D)-\LLb(x,D)\big)\bv(x)+\LLb(x,D)\bv(x)
$$
and using standard approach
(cf. \cite{CFF,DN}),  we obtain a representation of $\bv$ in terms of the
Newtonian-type  potentials 
$$
\bv(x)=\int_{B_r} \bs{\Gamma}(x_0;x-y) \LLb\bv(y) dy +
\int_{B_r} \bs{\Gamma}(x_0;x-y)\big(\LLb(x_0,D)-\LLb(y,D)\big) \bv(y) dy.
$$
Taking the $2b$-order  derivatives and then unfreezing the coefficients by
putting $x_0=x,$ we get
\begin{align}\label{10}
D^\alpha \bv(x)= &\  \underbrace{p.v.
\int_{B_r}D^\alpha\bs{\Gamma}(x;x-y)\LLb\bv(y)dy}_{\ds=:\K_\alpha (\LLb\bv)}\\
\nonumber
 &\ + \sum_{|\alpha'|=2b}  \underbrace{p.v.
\int_{B_r}D^{\alpha}\bs{\Gamma}(x;x-y)
	\big(\A_{\alpha'}(x)- \A_{\alpha'}(y)\big) D^{\alpha'}
	\bv(y)dy}_{\ds =:\CM_\alpha[\A_{\alpha'}, D^{\alpha'} \bv]}\\
\nonumber
&\ +   \int_{\SF^{n-1}}D^{\beta^s}\bs{\Gamma}(x;y)\nu_s
d\sigma_{y}\,\LLb\bv(x)\qquad \forall\ \alpha\colon |\alpha|=2b
\end{align}
where the derivatives  $D^\alpha\bs{\Gamma}(\cdot;\cdot)$ are taken with respect
to the second variable,  the multiindices $\beta^s$ are such that
$\beta^s:=(\alpha_1,\ldots,
\alpha_{s-1},\alpha_{s}-1,\alpha_{s+1},\ldots,\alpha_n),$ $|\beta^s|=2b-1$
and $\nu_s$ is the $s$-th component of the outer normal to $\SF^{n-1}.$

 Noting that each entry of the matrix $D^\alpha\bs{\Gamma}(x;y),$ $|\alpha|=2b,$ is a {\it
Calder\'on-Zygmund kernel\/} (cf. \cite{CFF,CFL}), we have
\begin{lem}\label{pr1}
Let $|\alpha|=|\alpha'|=2b$ and  $\A_\alpha\in L^\infty(\Omega).$
For each $p\in(1,\infty)$ and each $\lambda\in(0,n)$ there is a constant
$C$ depending on $n,$ $m,$ $b,$ $\delta,$ $\|\A_\alpha\|_{\infty;\Omega},$
$p$ and $\lambda$ such that
\begin{equation}\label{eqKfMor}
\|\K_\alpha \mathbf{f}\|_{p,\lambda;\Omega}\leq C
\|\mathbf{f}\|_{p,\lambda;\Omega},\qquad
 \|\CM_\alpha [\A_{\alpha'}, \mathbf{f} ]\|_{p,\lambda;\Omega}\leq C
\|\A_{\alpha'}\|_{*;\Omega}
 \|\mathbf{f}\|_{p,\lambda;\Omega}
 \end{equation}
for all $ \bff\in L^{p,\lambda}(\Omega).$
 Moreover, let $\A_\alpha\in VMO(\Omega)\cap L^\infty(\Omega)$ with
$VMO$-modulus $\eta_{\A_\alpha}.$
Then  for each $\varepsilon>0$ there exists
$r_0=r_0(\varepsilon,\eta_{\A_\alpha})$ such that if $r<r_0$ we have
\begin{equation}\label{3.11}
\|\CM_\alpha  [\A_{\alpha'}, \mathbf{f} ]\|_{p,\lambda;B_r}\leq C \varepsilon
\|\mathbf{f}\|_{p,\lambda;B_r}
\end{equation}
for all $ B_r\Subset \Omega$ and all $\bff\in
L^{p,\lambda}(B_r).$
\end{lem}
Lemma~\ref{pr1} is proved in \cite[Theorem~2.1, Corollary~2.7]{PS1}
in the general case of Calder\'on-Zygmund's kernels of mixed homogeneity.

\section{ Proof of Theorem~\ref{th3}}
Fix an arbitrary  $x_0\in \text{supp\,}\bu$ and let $B_r:=\{x\in\R^n:
|x-x_0|<r\}.$
The operators  $\K_\alpha$ and $\CM_\alpha [\A_{\alpha'},\cdot]$ are bounded
from $L^p$
into itself (cf.  \cite{CFL}) and therefore the representation \eqref{10} still
holds true (a.e. in $\Omega$)
for any function  $\bv\in
 W^{2b,p}_{0}(\Omega):=\mathrm{closure}_{W^{2b,p}(\Omega)}C_0^\infty(\Omega)$
and moreover
 for $\bv\in W^{2b,p}_{0}(\Omega)\cap W^{2b,p,\lambda}_{\rm loc}(\Omega)$ as
well. Let $\text{supp\,}\bv\subset B_r.$
In  view of \eqref{10}, \eqref{eqKfMor} and \eqref{3.11}, for each
$\varepsilon>0$ there exists
$r_0(\varepsilon,\eta_{\A_\alpha})$ such that for $r<r_0$ one has
$$
 \|D^{2b}\bv\|_{p,\lambda;B_r}\leq C(\|\LLb\bv\|_{p,\lambda;B_r}
+\varepsilon\|D^{2b}\bv\|_{p,\lambda;B_r}) $$
whence, choosing $\varepsilon$ small enough  we obtain
\begin{equation}\label{eq10}
\|D^{2b}\bv\|_{p,\lambda;B_r}\leq C \|\LLb\bv\|_{p,\lambda;B_r}.
\end{equation}
Let  $\theta\in(0,1),$ $\theta'=\theta(3-\theta)/2>\theta$ and define a cut-off
function  $\varphi(x)\in C_0^\infty(B_r)$ such that $\varphi(x)=1$ for $x\in
B_{\theta r}$ whereas
$\varphi(x)=0$	for $x\not\in B_{\theta' r}.$ It is clear that $|D^s\varphi|\leq
 C(s)[\theta(1-\theta)r]^{-s}$ for any $1\leq s\leq 2b$ because of
$\theta'-\theta=\theta(1-\theta)/2.$
Applying  \eqref{eq10} to $\bv(x):=\varphi(x)\bu(x)\in W^{2b,p}_{0}(B_r)\cap
W^{2b,p,\lambda}(B_r),$ we get
\begin{align*}
\|D^{2b}\bu\|_{p,\lambda;B_{\theta r}}\leq&\
   \|D^{2b}\bv\|_{p,\lambda;B_{\theta' r}}\leq C\|\LLb\bv\|_{p,\lambda;B_{\theta'
r}}\\
\leq&\	 C\left(\|\bff\|_{p,\lambda;B_{\theta'r}}+
\sum_{s=1}^{2b-1}\frac{\|D^{2b-s}
\bu\|_{p,\lambda;B_{\theta'r}}}{[\theta(1-\theta)r]^s}+
	 \frac{\|\bu\|_{p,\lambda;B_{\theta'r}}}{[\theta(1-\theta)r]^{2b}}
\right). \end{align*}
Hence, the choice of $\theta'$ and $\theta(1-\theta)\leq 2\theta'(1-\theta')$ imply
\begin{align*}
[\theta(1-\theta)r]^{2b}\|D^{2b}\bu\|_{p,\lambda;B_{\theta r}}
\leq&\	  C\Big([\theta'(1-\theta')r]^{2b}\|\bff\|_{p,\lambda;B_{\theta' r}}\\
& +   \sum_{s=1}^{2b-1}  [\theta'(1-\theta')r]^s
\|D^{s}\bu\|_{p,\lambda;B_{\theta'r}}+ \|\bu\|_{p,\lambda;B_{\theta'r}}\Big)
\end{align*}
Setting $\Theta_s$ for the weighted Morrey seminorms
$\sup_{0<\theta<1}[\theta(1-\theta)r]^s\|D^s\bu\|_{p,\lambda;B_{\theta r}},$
$s\in\{0,\ldots,2b\},$
the last inequality rewrites as
\begin{equation}\label{eq11}
\Theta_{2b}\leq C\left(r^{2b}\|\bff\|_{p,\lambda;B_r}+\sum_{s=1}^{2b-1}
\Theta_s+\Theta_0 \right).
\end{equation}
To manage the seminorms $\Theta_s$ with $1\leq s\leq 2b-1,$ we use the next
interpolation inequality which
follows from \cite[(5.6)]{Sl} as in \cite[Proposition~3.2]{PS1}:
\begin{prop}
There  exists a constant $C=C(n,m,b,p,\lambda,s)$ independent of $r$ and such
that $$
\Theta_s\leq \varepsilon
\Theta_{2b}+\frac{C}{\varepsilon^{s/(2b-s)}}\Theta_0\quad
\forall\  \varepsilon\in(0,2).
$$
\end{prop}
Therefore,  interpolating the intermediate seminorms in \eqref{eq11} and fixing
$\theta=1/2,$ we obtain the
following {\it Caccioppoli-type\/} inequality
\begin{equation}\label{eq12}
\|D^{2b}\bu\|_{p,\lambda;B_{r/2}}\leq C\big(\|\bff\|_{p,\lambda;B_r}+
r^{-2b}\|\bu\|_{p,\lambda;B_r}\big).
\end{equation}
The estimate \eqref{9a} follows from \eqref{eq12} by a finite covering of
$\Omega'$ with balls
$B_{r/2},$  $r<\text{dist\,}(\Omega',\partial\Omega'').$

\section{Proof of Corollary~\ref{cr4}}
To begin with, let $s=2b-1$ and $\Omega'_r=B_r\cap \Omega'$  with $2r<{\rm
dist\,}(\Omega',\partial\Omega'').$
 Direct calculations based on the classical Poincar\'e inequality  lead to
\begin{equation}\label{camp}
\frac{1}{r^{p+\lambda} }
\int_{\Omega'_r}|D^{2b-1}\bu(x)-(D^{2b-1}\bu)_{\Omega'_r}|^p
dx \leq C(n,p,m)  \|D^{2b}\bu\|^p_{p,\lambda;\tl \Omega'}
\end{equation}
 with suitable $\tl\Omega'$ such that $\Omega'\Subset\tl\Omega'\Subset\Omega''$
and $(\bff)_{\Omega'_r}$ standing
for the integral average of $\bff$ over $\Omega'_r.$
Taking the supremum with respect to $r$ we get that $D^{2b-1}\bu$ belongs to
the Campanato space
$\LL^{p,p+\lambda}(\Omega').$ Now, employing
the  embedding properties of Campanato spaces into Morrey, $BMO$ and H\"older
ones
(cf.  \cite[Theorem~2.1]{Cm}) and \eqref{9a}, we obtain as follows. If
$p+\lambda<n$ then
 $D^{2b-1}\bu\in L^{p,p+\lambda}(\Omega'),$ $p+\lambda=n$ implies
$D^{2b-1}\bu\in BMO(\Omega'),$
while $D^{2b-1}\bu\in C^{0,\sigma_{2b-1}}(\Omega')$ with
$\sigma_{2b-1}=1-(n-\lambda)/p$ when $p+\lambda>n.$

To prove the statement for any $s,$ we run induction
for  decreasing values of $s$ until $s_0.$ Thus, suppose $D^{s+1} \bu$ satisfies
the statement of Corollary~$\ref{cr4}.$ Then
\begin{align*}
\frac{1}{r^{(2b-s)p+\lambda}}
\int_{\Omega'_r} |D^{s}\bu(x)-(D^{s}\bu)_{\Omega'_r}|^p dx
\leq&\  \|D^{s+1}\bu\|^p_{p,(2b-s-1)p+\lambda;\tl \Omega'}\\
\leq&\ 
C\left( \|\bff\|^p_{p,\lambda;\Omega''}+\|\bu\|^p_{p,\lambda;\Omega''}\right)
\end{align*}
whence	$D^s\bu\in \LL^{p,(2b-s)p+\lambda}(\Omega')$ and the conclusions follow
as above.

Similarly, if $s_0\geq 1$ and $p\in \left(1,\frac{n-\lambda}{2b-s_0}\right)$
we have  $D^{s_0-1}\bu\in \LL^{p,(2b-s_0+1)p+\lambda}(\Omega').$ It follows from
the definition of $s_0$ that $n<(2b-s_0+1)p+\lambda< n+p,$
whence $D^{s_0-1}\bu\in C^{0,2b-s_0+1-\frac{n-\lambda}{p}}(\Omega').$


\begin{thebibliography}{22}
\bibitem{Cm}
{S. Campanato},
{\em Sistemi ellittici in forma divergenza. Regolarit\`a
all'interno,} Pubblicazioni della Classe di Scienze: Quaderni,
Scuola Norm. Sup., Pisa, 1980.
\bibitem{Cm1}  S. Campanato,
Propriet\'a di H\"olderianit\`a di alcune classi
di funzioni.
{\em Ann.\ Scuola\ Norm.\ Sup.\ Pisa }{\textbf 17}, (1963) 175-188.
\bibitem{CFF}
{F. Chiarenza, M.~Franciosi \and  M.~Frasca},
$L^p$-estimates for linear
elliptic systems with discontinuous coefficients,
{\em Rend.  Accad.  Naz.  Lincei,  Mat.  Appl.,}~{\bf 5}~(1994)~27--32.
\bibitem{CFL}
{F. Chiarenza,	M. Frasca \and P.~Longo},
Interior $W^{2,\,p}$ estimates for non divergence elliptic
equations with discontinuous coefficients,
{\em Ric.  Mat.,} {\bf 60} (1991) 149--168.
\bibitem{DN}
{A. Douglis \and L. Nirenberg},
Interior estimates for elliptic systems of partial differential equations,
{\em Commun.  Pure  Appl.  Math.,} {\bf 8} (1955) 503--538.
\bibitem{J}
{F.~John}, {\em Partial Differential Equations,}
Appl. Math. Sci., Vol. 1, Springer-Verlag, Berlin, 1991.
\bibitem{JN}
{F.~John  \and	L.~Nirenberg},
On functions of bounded mean oscillation,
{\em Commun.  Pure  Appl.  Math.,} {\bf 14} (1961) 415--426.
\bibitem{MPS}
{A. Maugeri, D.K. Palagachev \and L.G. Softova},
{\em Elliptic  and  Parabolic  Equations  with	Discontinuous  Coefficients,}
Wiley-VCH, Berlin, 2000.
\bibitem{M} C.B.  Morrey,  Jr.,
{\em Multiple Integrals in the Calculus of
Variations.} Springer-Verlag, Berlin, 1966.
\bibitem{PS1}
{D. Palagachev \and  L. Softova},
Singular integral operators, Morrey spaces and fine
regularity of solutions to PDE's, {\em Potential  Anal.,} {\bf 20}
(2004) 237--263.
\bibitem{S}
{D. Sarason},
Functions of vanishing mean oscillation,
{\em  Trans.  Amer.  Math.  Soc.,} {\bf 207} (1975) 391--405.
\bibitem{Sl}
{V. A. Solonnikov},
On the boundary value problems for linear parabolic systems of differential
equations of general form,
{\em Proc.  Steklov  Inst.  Math.,} {\bf 83} (1965);
English translation: O. A. Ladyzhenskaya~(Ed.)
{\em Boundary  Value  Problems	of  Mathematical  Physics  III,}
Amer. Math. Soc., Providence, R.I., 1967.
\end{thebibliography}
\end{document}